\documentclass[12pt]{article}
\usepackage{mathrsfs}
\usepackage{amsthm}
\usepackage{amssymb}
\usepackage{latexsym}
\usepackage{amsmath,amsfonts}
\numberwithin{equation}{section}
\usepackage{cases}
\usepackage{bm}
\usepackage{indentfirst}
\usepackage{xcolor}
\usepackage{ifpdf}
\usepackage{graphicx}
\usepackage{epstopdf}
\usepackage{epsfig}
\usepackage{psfrag}
\usepackage{enumitem}
\usepackage{verbatim}
\usepackage{color}
\usepackage{subfigure}

\usepackage[
pdfauthor={Lu},
pdftitle={matching covered},
pdfstartview=XYZ,
bookmarks=true,
colorlinks=true,
linkcolor=blue,
urlcolor=blue,
citecolor=blue,
bookmarks=true,
linktocpage=true,
hyperindex=true
]{hyperref}

\usepackage[textwidth=16cm, textheight=22cm, includehead, includefoot]{geometry}

\title{The maximum number of edges in minimal $k$-extendable bipartite graphs \footnote{E-mail address: xiaolinghe99@163.com.} }
\author{ Xiaoling He\\
\small {School of Mathematics and Statistics, Lanzhou University, Lanzhou, China}
}
\date{}

\newtheorem{lem}{Lemma}[section]
\newtheorem{thm}[lem]{Theorem}

\newtheorem{pro}[lem]{Proposition}
\newtheorem{conj}[lem]{Conjecture}

\newtheorem*{remark}{Remark}

\newtheorem*{claim}{Claim}

\begin{document}
\bibliographystyle{plain}
\newcommand{\udots}{\mathinner{\mskip1mu\raise1pt\vbox{\kern7pt\hbox{.}}
\mskip2mu\raise4pt\hbox{.}\mskip2mu\raise7pt\hbox{.}\mskip1mu}}
\maketitle

\begin{abstract}
    A connected graph of order at least $2k+2$ is {\em $k$-extendable} if it has a matching of size $k$ and every such matching extends to a perfect matching. 
    In particular, a 1-extendable graph is also called a {\em matching covered graph}.
    A $k$-extendable graph $G$ is {\em minimal} if $G-e$ is not $k$-extendable for every $e\in E(G)$.
    Lov\'asz and Plummer [J. Combin. Theory Ser. B 23 (1977), 127--138] proved that every minimal matching covered bipartite graph $G$ with $|V(G)|\ge6$ satisfies $|E(G)|\le(3|V(G)|-6)/2$, and this bound is sharp for every even order at least 6. 
    Fabres, Kothari and Carvalho [J. Graph Theory 96 (2021), 490--509] proved that every minimal 2-extendable bipartite graph $G$ with $|V(G)|\ge12$ satisfies $|E(G)|\le(5|V(G)|-20)/2$. 
    Very recently,  Mallik, Diwan and Kothari [J. Combin. 17 (2026), 559--606] characterized minimal matching covered bipartite graphs $G$ with $|E(G)|=(3|V(G)|-6)/2$ in terms of special trees. 
    They also  conjectured that, for every integer $k\ge1$, there exists an integer $N_k\le4k^2+2k$ such that every minimal $k$-extendable bipartite graph $G$ with $|V(G)|\ge N_k$ satisfies $|E(G)|\le(2k+1)(|V(G)|-2k)/2$. 
    In this paper, we confirm this conjecture for $k\ge3$ and prove that the upper bound is sharp for every even order at least $6k$ ($k\ge3$). 
\end{abstract}
\vskip 0.2in \baselineskip 0.1in

\section{Introduction}
Graphs considered in this paper may have multiple edges, but no loops. 
We follow \cite{BM08} for undefined notation and terminology.
For a graph $G$, let $V(G)$ and $E(G)$ denote its vertex set and edge set, respectively. 
Denote by $G[A,B]$ the bipartite graph with bipartition $(A,B)$. 
Let $G_1$ and $G_2$ be subgraphs of $G$. The union $G_1\cup G_2$ of $G_1$ and $G_2$ is the subgraph with vertex set $V(G_1)\cup V(G_2)$ and edge set $E(G_1)\cup E(G_2)$. 

A {\em perfect matching} $M$ of a graph $G$ is a set of edges such that each vertex of $G$ is incident with exactly one edge of $M$.
For an integer $k\ge1$, a connected graph of order at least $2k+2$ is {\em $k$-extendable} if it has a matching of size $k$ and every such matching extends to a perfect matching.
A $k$-extendable graph $G$ is {\em minimal} if $G-e$ is not $k$-extendable for every $e\in E(G)$.

The concept of $k$-extendable graphs was introduced by Plummer \cite{Plummer1980}.
Anunchuen and Caccetta \cite{AC1994,AC1997} studied minimal $k$-extendable graphs.
A 1-extendable graph is also called a {\em matching covered} graph. 
Lov\'asz and Plummer \cite{LP77} proved the following result for minimal matching covered bipartite graphs.\footnote{Lov\'asz and Plummer used the terminology ``minimal elementary bipartite graph'' in \cite{LP77}.}

\begin{thm}[\cite{LP77}]\label{thm:minimal-bi-MC}
Every minimal matching covered bipartite graph $G$ with $|V(G)|\ge6$ satisfies $|E(G)|\le(3|V(G)|-6)/2$, and this bound is sharp for every even order at least $6$.
\end{thm}

Very recently, Mallik, Diwan and Kothari \cite{MDK2026} provided a complete characterization of the minimal matching covered bipartite graphs attaining the upper bound in Theorem \ref{thm:minimal-bi-MC}, and proposed the following conjecture.

\begin{conj}[{\cite[Conjecture~8.11]{MDK2026}}]\label{conj:MDK2026}
For every integer $k\ge1$, there exists an integer $N_k\le4k^2+2k$ such that every minimal $k$-extendable bipartite graph $G$ with $|V(G)|\ge N_k$ satisfies $|E(G)|\le(2k+1)(|V(G)|-2k)/2$.
\end{conj}

Theorem \ref{thm:minimal-bi-MC} confirms Conjecture \ref{conj:MDK2026} for $k=1$.
A 2-extendable bipartite graph is also called a {\em brace}.
Fabres, Kothari and Carvalho \cite{FKC2020} proved the following result, which confirms Conjecture \ref{conj:MDK2026} for $k=2$.

\begin{thm}[\cite{FKC2020}]\label{minimal_brace}
Every minimal brace $G$ with $|V(G)|\ge12$ satisfies $|E(G)|\le(5|V(G)|-20)/2$, and this bound is sharp for every even order at least $12$.
\end{thm}

In this paper, we prove the following sharp upper bound on the number of edges in minimal $k$-extendable bipartite graphs for $k\ge3$.

\begin{thm}\label{thm:main}
Let $G$ be a minimal $k$-extendable bipartite graph with $k\ge3$ and $|V(G)|\ge6k$.
Then $|E(G)|\le(2k+1)(|V(G)|-2k)/2$, and this bound is sharp for every even order at least $6k$.
\end{thm}

Together with Theorems \ref{thm:minimal-bi-MC} and \ref{minimal_brace}, Theorem \ref{thm:main} confirms Conjecture \ref{conj:MDK2026} with $N_k=6k$ for every integer $k\ge1$.

\begin{figure}[!h]
\centering
\includegraphics[height=3.5cm]{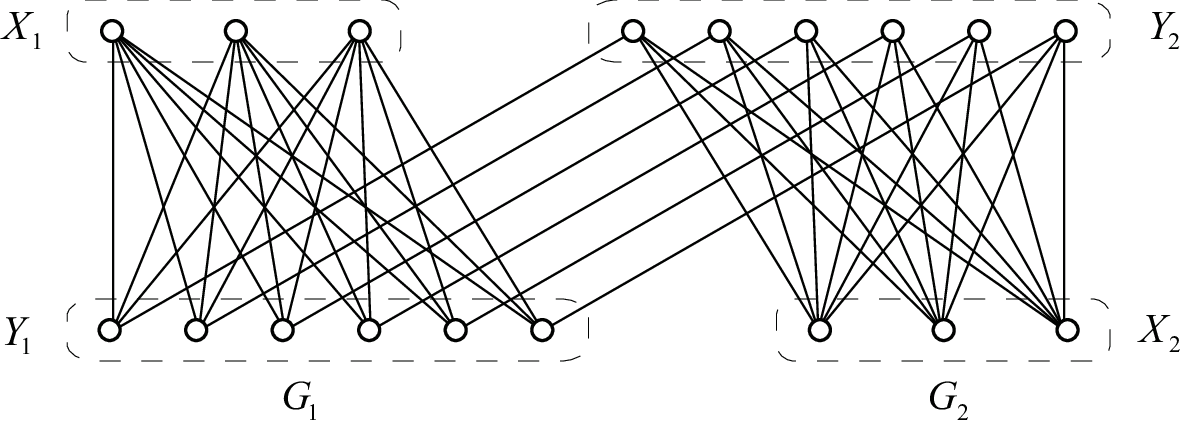}
\caption{The graph $G_{3,18}$.}
\label{fig:G_kn}
\end{figure}

For any integer $k\ge3$ and any even integer $n\ge6k$, let $G_1[X_1,Y_1]$ and $G_2[X_2,Y_2]$ be two disjoint copies of $K_{k,n/2-k}$ such that $|X_1|=|X_2|=k$ and $|Y_1|=|Y_2|=n/2-k$. 
Let $G_{k,n}$ be the graph obtained from $G_1\cup G_2$ by adding a perfect matching between $Y_1$ and $Y_2$ (see Figure \ref{fig:G_kn} for $G_{3,18}$);
equivalently, $G_{k,n}$ is obtained from the star
$K_{1,n/2-k}$ by isomorphic $k$-leaf matching
(see \cite[Definition 8.18]{MDK2026}).
Then $G_{k,n}$ is bipartite with $X_1\cup Y_2$ and $Y_1\cup X_2$ as its two partite sets.   
The following result, proved by Mallik, Diwan and Kothari, implies that the upper bound in Theorem \ref{thm:main} is sharp.

\begin{lem}[{\cite[Theorem 8.25 (ii)]{MDK2026}}]\label{lem:sharp}
For any integer $k\ge3$, any graph $G$ obtained from a star $K_{1,p}$ with $p\ge2k$ by isomorphic $k$-leaf matching is a minimal $k$-extendable bipartite graph and satisfies $|E(G)|=(2k+1)(|V(G)|-2k)/2$.
\end{lem} 
Since $n\ge6k$, we have $n/2-k\ge2k$ and hence,  $G_{k,n}$ is a minimal $k$-extendable bipartite graph of order $n$ with $|E(G_{k,n})|=(2k+1)(n-2k)/2$ by Lemma \ref{lem:sharp}. 

\begin{remark}
{\rm  
The condition $|V(G)|\ge6k$ in Theorem \ref{thm:main} cannot be weakened to $|V(G)|\ge6k-2$. 
For each integer $k\ge3$ and each even integer $n$ with $2k+2\le n\le6k-2$, we define the graph $H_{k,n}$ in Section \ref{sec:6k}, and prove that it is a minimal $k$-extendable bipartite graph of order $n$ with $|E(H_{k,n})|>(2k+1)(n-2k)/2$ (see Proposition \ref{pro:6k}).
}
\end{remark}

In Section \ref{sec:Preliminaries}, we present some notation and preliminary results. In Sections \ref{sec:proof_of_main_result} and \ref{sec:6k}, we prove Theorem \ref{thm:main}, and prove that $|V(G)|\ge 6k$ cannot be weakened to $|V(G)|\ge 6k-2$, respectively.

\section{Preliminaries}\label{sec:Preliminaries}
A set of vertices $S$ of $G$ is a {\em vertex cover} of $G$ if each edge of $G$ has at least one end vertex in $S$. 
Let $\nu(G)$ denote the maximum size of a matching in $G$. 
 
For $\emptyset\neq X\subset V(G)$, let $N_G(X)$ denote the set of all vertices in $V(G)\setminus X$ which are adjacent to at least one vertex in $X$. If $X=\{u\}$, we write $N_G(u)$, or simply $N(u)$, for $N_G(X)$.  

For $X\subseteq V(G)$, let $G[X]$ be the subgraph of $G$ induced by $X$. 
For disjoint subsets $X,Y\subseteq V(G)$, let $E_G[X,Y]$ denote the set of edges with one end vertex in $X$ and the other in $Y$. 
For $X\subseteq V(G)$, let $\partial_G(X):=E_G[X,V(G)\setminus X]$.

For a vertex $u\in V(G)$, let $d_G(u)$, or simply $d(u)$, denote the degree of $u$ in $G$, that is, the number of edges of $G$ incident with $u$.
A vertex of degree zero is called an {\em isolated vertex}. 
Let $\delta(G)$ denote the {\em minimum degree} of $G$.

For a minimal $k$-extendable bipartite graph $G:=G[A,B]$ with $|V(G)|\ge2k+2$, let
\begin{equation*}
\begin{aligned}
V_{k+1}:=\{u\in V(G):d_G(u)=k+1\}, ~
&V_{k+2}:=\{u\in V(G):d_G(u)\ge k+2\},\\ 
E_{k+1}:=E(G[V_{k+1}]),~ &E_{k+2}:=E(G[V_{k+2}]),\\ 
A_{k+1}:=A\cap V_{k+1},~ &B_{k+1}:=B\cap V_{k+1},\\ 
A_{k+2}:=A\cap V_{k+2},~ &B_{k+2}:=B\cap V_{k+2}.
\end{aligned}
\end{equation*}  
Since $G$ has a perfect matching, we have $|A|=|B|$.
Clearly, every minimal $k$-extendable graph is simple for any integer $k\ge 1$.

\begin{thm}[\cite{Plummer1980}]\label{thm:con_k-ext}
Let $k\ge1$ be an integer and let $G$ be a $k$-extendable graph. 
Then $G$ is $(k+1)$-connected. If $k\ge2$, then $G$ is also $(k-1)$-extendable.
\end{thm}

For any integer $k\ge1$, every $k$-extendable graph $G$ satisfies $\delta(G)\ge k+1$ by Theorem \ref{thm:con_k-ext}. 
By Theorem \ref{thm:con_k-ext}, every matching of size at most $k$ in a $k$-extendable graph extends to a perfect matching. 
Moreover, Lou \cite{Lou1999} proved the following result for minimal $k$-extendable bipartite graphs.

\begin{lem}[\cite{Lou1999}]\label{lem:min_deg}
Let $G$ be a minimal $k$-extendable bipartite graph. 
Then $\delta(G)=k+1$ and $G[V_{k+2}]$ is a forest.
Consequently, if $V_{k+2}\neq\emptyset$, then
\begin{equation}\label{eq:forest}
    |E_{k+2}|\le |V_{k+2}|-1.
\end{equation}
\end{lem}

\begin{thm}[\cite{Plummer1986}]\label{thm:char_k-ext-bi}
Let $k\ge1$ be an integer, and let $G:=G[A,B]$ be a connected bipartite graph with $|V(G)|\ge 2k+2$, where $|A|=|B|$. Then the following are equivalent.\\ 
{\rm (i)}  $G$ is $k$-extendable.\\ 
{\rm (ii)} For every nonempty $S\subseteq A$, either $N_G(S)=B$ or $|N_G(S)|\ge |S|+k$; equivalently, $|N_G(S)|\ge |S|+k$ for every $S\subseteq A$ with $1\le |S|\le |A|-k$.\\ 
{\rm (iii)} $G-A'-B'$ has a perfect matching for all $A'\subseteq A$ and $B'\subseteq B$ such that $|A'|=|B'|=k$. 
\end{thm}
 
\begin{lem}\label{lem:edge_count}
Let $G:=G[A,B]$ be a minimal $k$-extendable bipartite graph. Then\\ 
{\rm (i)} $|E(G)|=(k+1)|V_{k+1}|+|E_{k+2}|-|E_{k+1}|$, and \\ 
{\rm (ii)} $|E_{k+1}|\ge|A|-|V_{k+2}|+\min\{k,\nu(G[V_{k+2}])\}$. 
\end{lem}
\begin{proof} 
Since every vertex of $V_{k+1}$ has degree $k+1$, we have $|E(G)|=|E_{k+1}|+|\partial_G(V_{k+1})|+|E_{k+2}|=(k+1)|V_{k+1}|+|E_{k+2}|-|E_{k+1}|$ and so, (i) holds.

Let $M_2$ be a matching of $G[V_{k+2}]$ such that $|M_2|=\min\{k,\nu(G[V_{k+2}])\}$. 
Since $|M_2|\le k$, $M_2$ can be extended to a perfect matching $M$ of $G$ by Theorem \ref{thm:con_k-ext}.   
Note that $2|M\cap E_{k+2}|+|M\cap\partial_G(V_{k+2})|=|V_{k+2}|$. 
Since $|M|=|A|$ and $M_2\subseteq M\cap E_{k+2}$, we have $|E_{k+1}|\ge |M\cap E_{k+1}|=|A|-|V_{k+2}|+|M\cap E_{k+2}|\ge |A|-|V_{k+2}|+|M_2|$ and so, (ii) holds. 
\end{proof}

\section{Proof of Theorem \ref{thm:main}}\label{sec:proof_of_main_result}

\begin{lem}\label{lem:E2}
Let $G:=G[A,B]$ be a minimal $k$-extendable bipartite graph, and let $uv\in E_{k+2}$ with $u\in A$.
Then there exist $A_{uv}\subseteq A$ and $B_{uv}\subseteq B$ such that $1\le|A_{uv}|,|B_{uv}|\le|A|-k$ and\\ 
{\rm (i)} $E_G[A_{uv},B_{uv}]=\{uv\}$;\\
{\rm (ii)} $|A_{uv}|+|B_{uv}|=|A|-k+1$;\\
{\rm (iii)} $|N_G(A_{uv})|=|A_{uv}|+k$ and $|N_G(B_{uv})|=|B_{uv}|+k$. 
\end{lem}
\begin{proof}
Since $G$ is minimal $k$-extendable, $G-uv$ is not $k$-extendable. 
By Theorem \ref{thm:con_k-ext}, $G$ is 2-connected and so $G-uv$ is connected.
By Theorem \ref{thm:char_k-ext-bi} (ii), there exists a nonempty subset $A_{uv}\subseteq A$ with $|A_{uv}|\le|A|-k$ such that $|N_{G-uv}(A_{uv})|\le|A_{uv}|+k-1$. 
Since $G$ is $k$-extendable, we have $|A_{uv}|+k\le |N_G(A_{uv})|\le |N_{G-uv}(A_{uv})|+1\le |A_{uv}|+k$. 
Thus, $u\in A_{uv}$, $N_G(v)\cap A_{uv}=\{u\}$, and
\begin{equation}\label{eq:NXA}
|N_G(A_{uv})|=|A_{uv}|+k
~\mbox{and}~
|N_{G-uv}(A_{uv})|=|A_{uv}|+k-1.
\end{equation}
Let $B_{uv}:=B\setminus N_{G-uv}(A_{uv})$.
Then $v\in B_{uv}$ and $E_G[A_{uv},B_{uv}]=\{uv\}$. So (i) holds.
By \eqref{eq:NXA} and $|A|=|B|$, we have
$|A_{uv}|+|B_{uv}| =|A_{uv}|+|B|-|N_{G-uv}(A_{uv})| =|A|-k+1$ and so, (ii) follows.  

Since $1\le |B_{uv}|\le |B|-k$, we have $|N_G(B_{uv})|\ge |B_{uv}|+k$ by Theorem \ref{thm:char_k-ext-bi} (ii). 
By (i), $N_G(B_{uv})\subseteq (A\setminus A_{uv})\cup\{u\}$.
Hence, by (ii), $|N_G(B_{uv})|\le |A|-|A_{uv}|+1=|B_{uv}|+k$.  
Together with \eqref{eq:NXA}, (iii) follows. 
\end{proof}

We call the sets $A_{uv}$ and $B_{uv}$ in Lemma \ref{lem:E2} a pair of {\em end sets}\footnote{Lou used the terminology ``key set'' in \cite{Lou1999}.} of $uv$.

\begin{lem}\label{lem:adjacent_E2}
Let $G:=G[A,B]$ be a minimal $k$-extendable bipartite graph, and let $uv_1,uv_2\in E_{k+2}$ be two distinct edges incident with $u\in A$. 
For each $i\in\{1,2\}$, let $A_{uv_i},B_{uv_i}$ be a pair of end sets of $uv_i$.
Then $B_{uv_1}\cap B_{uv_2}=\emptyset$.
\end{lem}
\begin{proof}
Suppose to the contrary that $B_{uv_1}\cap B_{uv_2}\neq\emptyset$. 
By Lemma \ref{lem:E2} (i), we have $v_i\in B_{uv_i}\setminus B_{uv_{3-i}}$ and $N_G(u)\cap B_{uv_i}=\{v_i\}$ for $i\in\{1,2\}$.
Thus, $u\in N_G(B_{uv_1})\cap N_G(B_{uv_2})$, but
$u\notin N_G(B_{uv_1}\cap B_{uv_2})$.
Moreover, $N_G(u)\cap(B_{uv_1}\cup B_{uv_2})=\{v_1,v_2\}$.
Since $G$ is simple and $d_G(u)\ge k+2$, we have
\begin{equation}\label{eq:end-union}
|B_{uv_1}\cup B_{uv_2}|\le |B|-d_G(u)+2\le |B|-k.
\end{equation}  
Since $u\in N_G(B_{uv_1})\cap N_G(B_{uv_2})$ and
$u\notin N_G(B_{uv_1}\cap B_{uv_2})$, by \eqref{eq:end-union} and Theorem \ref{thm:char_k-ext-bi} (ii), we have
\begin{equation}\label{eq:end-neighbor}
\begin{aligned}
|N_G(B_{uv_1})\cap N_G(B_{uv_2})|
&\ge |N_G(B_{uv_1}\cap B_{uv_2})|+1\ge |B_{uv_1}\cap B_{uv_2}|+k+1,\\
|N_G(B_{uv_1}\cup B_{uv_2})|
&\ge |B_{uv_1}\cup B_{uv_2}|+k.
\end{aligned}
\end{equation} 
By Lemma \ref{lem:E2} (iii) and \eqref{eq:end-neighbor}, we have 
\begin{equation*}
\begin{aligned}
|B_{uv_1}|+|B_{uv_2}|+2k&=|N_G(B_{uv_1})|+|N_G(B_{uv_2})|~~~(\mbox{by Lemma \ref{lem:E2} (iii)})\\
&=|N_G(B_{uv_1}\cup B_{uv_2})|+|N_G(B_{uv_1})\cap N_G(B_{uv_2})|\\
&\ge |B_{uv_1}\cup B_{uv_2}|+|B_{uv_1}\cap B_{uv_2}|+2k+1~~~(\mbox{by \eqref{eq:end-neighbor}})\\
&=|B_{uv_1}|+|B_{uv_2}|+2k+1.
\end{aligned}
\end{equation*} 
It is a contradiction. Therefore, the lemma follows.
\end{proof}

\begin{figure}[!h]
\centering
\includegraphics[height=3.5cm]{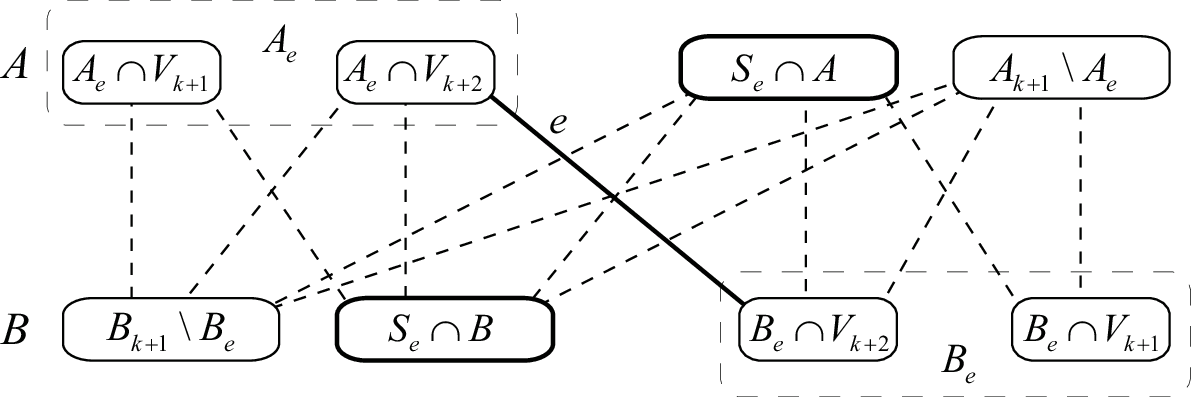}
\caption{The sets $A_e$, $B_e$, and $S_e$ for $e\in E_{k+2}$.}
\label{fig:Se}
\end{figure}

By symmetry, if $u_1v,u_2v\in E_{k+2}$ are two distinct edges such that $v\in B$, then $A_{u_1v}\cap A_{u_2v}=\emptyset$. 
For every $e\in E_{k+2}$, choose a pair of end sets $A_{e},B_{e}$, and let $S_e:=V_{k+2}\setminus(A_e\cup B_e)$ (see Figure \ref{fig:Se}).

\begin{lem}\label{lem:lowE1}
Let $G:=G[A,B]$ be a minimal $k$-extendable bipartite graph with $|V(G)|\ge 2k+2$.
Then
\begin{equation}\label{eq:Se}
\begin{aligned}
2|E_{k+1}|\ge{}&
(|A|-k+1-|V_{k+2}|)|E_{k+2}|+(2k+2-|V_{k+2}|)|A|
+2|A_{k+2}||B_{k+2}|\\
&-(k+1)|V_{k+2}|+\sum_{e\in E_{k+2}}|S_e|.
\end{aligned}
\end{equation}
Moreover,\\
{\rm (i)}$\sum_{e\in E_{k+2}}|S_e|\ge\max\{(|E_{k+2}|-1)\nu(G[V_{k+2}]),\sum_{u\in V_{k+2}}
d_{G[V_{k+2}]}(u)(d_{G[V_{k+2}]}(u)-1)\}$;\\
{\rm (ii)} if $|E_{k+2}|\ge5$, then
$\sum_{e\in E_{k+2}}|S_e|\ge3|E_{k+2}|-3$.
\end{lem}
\begin{proof} 
We first prove \eqref{eq:Se}.
Let $v\in B_{k+2}$.
By Lemma \ref{lem:adjacent_E2}, the sets
$A_{uv}\cap V_{k+1}$, over all $uv\in E_{k+2}$
incident with $v$, are pairwise disjoint.
For every $uv\in E_{k+2}$ incident with $v$, we have $N_G(v)\cap A_{uv}=\{u\}$ by Lemma \ref{lem:E2} (i).
Since $u\in V_{k+2}$, we have $N_G(v)\cap(A_{uv}\cap V_{k+1})=\emptyset$ and so,  
\begin{equation}\label{eq:Nv}
|N_G(v)\cap V_{k+1}|\le|A_{k+1}|-\sum_{\substack{u\in A\\uv\in E_{k+2}}}|A_{uv}\cap V_{k+1}|.
\end{equation}
Since every vertex of $A_{k+1}$ has degree $k+1$, 
summing \eqref{eq:Nv} over $v\in B_{k+2}$, we have 
\begin{equation}\label{eq:low-1}
\begin{aligned}
|E_{k+1}|
&=(k+1)|A_{k+1}|
-\sum_{v\in B_{k+2}}|N_G(v)\cap V_{k+1}|\\
&\ge(k+1-|B_{k+2}|)|A_{k+1}|
+\sum_{e\in E_{k+2}}|A_e\cap V_{k+1}|.
\end{aligned}
\end{equation}
By symmetry, we have
\begin{equation}\label{eq:low-2}
|E_{k+1}|\ge(k+1-|A_{k+2}|)|B_{k+1}|+\sum_{e\in E_{k+2}}|B_e\cap V_{k+1}|.
\end{equation}
For every $e\in E_{k+2}$, by Lemma \ref{lem:E2} (ii) and the definition of $S_e$, we have
\begin{equation}\label{eq:k1}
\begin{aligned}
|A_e\cap V_{k+1}|+|B_e\cap V_{k+1}|
&=|A_e|+|B_e|-|(A_e\cup B_e)\cap V_{k+2}|\\
&=|A|-k+1-|V_{k+2}|+|S_e|.
\end{aligned}
\end{equation} 
Adding \eqref{eq:low-1} and \eqref{eq:low-2}, and using \eqref{eq:k1} together with $|A_{k+1}|=|A|-|A_{k+2}|$ and $|B_{k+1}|=|A|-|B_{k+2}|$, we obtain \eqref{eq:Se}. 

(i)
For every $e\in E_{k+2}$, let $w\in A$ and $w'\in B$ be the end vertices of $e$.
By Lemma \ref{lem:E2} (i), we have $E_G[A_e,B_e]=\{e\}$.
Hence, $w\in A_e$, $w'\in B_e$ and $w,w'\notin S_e$.
Since $G$ is bipartite, every edge of $G[V_{k+2}]$ with both end vertices in $A_e\cup B_e$ belongs to
$E_G[A_e,B_e]$. 
Thus, every edge of $G[V_{k+2}]-e$ has an end vertex in $S_e$ and so, $S_e$ is a vertex cover of $G[V_{k+2}]-e$. 

Since $w,w'\notin S_e$ and $S_e$ covers every edge of $G[V_{k+2}]-e$, we have 
$N_{G[V_{k+2}]}(w)\setminus\{w'\}\subseteq S_e$ and
$N_{G[V_{k+2}]}(w')\setminus\{w\}\subseteq S_e$.
Moreover,
$N_{G[V_{k+2}]}(w)\setminus\{w'\}\subseteq B$ and
$N_{G[V_{k+2}]}(w')\setminus\{w\}\subseteq A$ are disjoint. 
Thus, 
\begin{equation}\label{eq:Se_deg}
|S_e|\ge d_{G[V_{k+2}]}(w)+d_{G[V_{k+2}]}(w')-2.
\end{equation}
Summing \eqref{eq:Se_deg} over $e\in E_{k+2}$, we have
\begin{equation}\label{eq:sum_Se}
\begin{aligned}
\sum_{e\in E_{k+2}}|S_e|&\ge \sum_{e=ww'\in E_{k+2}}(d_{G[V_{k+2}]}(w)+d_{G[V_{k+2}]}(w')-2)\\
&=\sum_{u\in V_{k+2}}d_{G[V_{k+2}]}(u)(d_{G[V_{k+2}]}(u)-1).
\end{aligned}
\end{equation}

Let $M$ be a maximum matching of $G[V_{k+2}]$.
For every $e\in E_{k+2}\setminus M$, $S_e$
covers every edge of $M$.
For every $f\in M$, $S_f$ covers every edge of
$M\setminus\{f\}$. 
Since $M$ is a matching,  we have
\begin{equation}\label{eq:Se_matching}
|S_e|\ge
\begin{cases}
|M|,& e\in E_{k+2}\setminus M,\\
|M|-1,& e\in M.
\end{cases}
\end{equation}
Summing \eqref{eq:Se_matching} over $e\in E_{k+2}$, we have
\begin{equation}\label{eq:sum_Se_matching}
\begin{aligned}
\sum_{e\in E_{k+2}}|S_e|
&\ge(|E_{k+2}|-|M|)|M|+|M|(|M|-1)\\
&=(|E_{k+2}|-1)\nu(G[V_{k+2}]).
\end{aligned}
\end{equation}
By \eqref{eq:sum_Se} and \eqref{eq:sum_Se_matching},
{\rm (i)} holds.

{\rm (ii)}
Assume that $|E_{k+2}|\ge5$.
If $\nu(G[V_{k+2}])\ge3$, then
$\sum_{e\in E_{k+2}}|S_e|
\ge(|E_{k+2}|-1)\nu(G[V_{k+2}])
\ge3|E_{k+2}|-3$ by \eqref{eq:sum_Se_matching}.
If $\nu(G[V_{k+2}])=1$, then, since $G[V_{k+2}]$
is bipartite, $G[V_{k+2}]$ consists of a star with
$|E_{k+2}|$ edges and possibly some isolated vertices.
By \eqref{eq:sum_Se} and $|E_{k+2}|\ge5$, we have
$\sum_{e\in E_{k+2}}|S_e|
\ge|E_{k+2}|(|E_{k+2}|-1)\ge3|E_{k+2}|-3$.

Now assume that $\nu(G[V_{k+2}])=2$.
By \cite[Theorem 1.1.1]{LP1986}, there exists a vertex
cover of $G[V_{k+2}]$ of size $2$. 
If $G[V_{k+2}]$ has a vertex cover $\{s,t\}$ with $st\in E_{k+2}$, then every edge of $G[V_{k+2}]$ is incident with $s$ or $t$, and $st$ is the only edge incident with both $s$ and $t$. Hence, 
\begin{equation}\label{eq:deg_sum}
d_{G[V_{k+2}]}(s)+d_{G[V_{k+2}]}(t)=|E_{k+2}|+1.
\end{equation}
By \eqref{eq:sum_Se}, \eqref{eq:deg_sum} and $|E_{k+2}|\ge5$,
we have $\sum_{e\in E_{k+2}}|S_e|\ge d_{G[V_{k+2}]}(s)(d_{G[V_{k+2}]}(s)-1) +d_{G[V_{k+2}]}(t)(d_{G[V_{k+2}]}(t)-1)\ge \frac{|E_{k+2}|^2-1}{2}\ge 3|E_{k+2}|-3.$

Now assume that $G[V_{k+2}]$ has no vertex cover
consisting of two adjacent vertices. 
Let $\{s,t\}$ be a vertex cover of $G[V_{k+2}]$. 
Then $st\notin E_{k+2}$. 
Since $\nu(G[V_{k+2}])=2$, we have $d_{G[V_{k+2}]}(s)\ge1$ and $d_{G[V_{k+2}]}(t)\ge1$. 
By Lemma \ref{lem:min_deg}, $G[V_{k+2}]$ is a forest. 
So $|N_{G[V_{k+2}]}(s)\cap N_{G[V_{k+2}]}(t)|\le 1$. 
Suppose that $N_{G[V_{k+2}]}(t)\subseteq N_{G[V_{k+2}]}(s)$. 
Then $t$ has exactly one neighbor $z$ in $G[V_{k+2}]$ as $d_{G[V_{k+2}]}(t)\ge1$. 
Since $\{s,t\}$ is a vertex cover of $G[V_{k+2}]$, every edge of $E_{k+2}\setminus\{tz\}$ is incident with $s$. 
So $\{z,s\}$ is also a vertex cover of $G[V_{k+2}]$ such that $zs\in E_{k+2}$, contradicting the assumption that $G[V_{k+2}]$ has no vertex cover
consisting of two adjacent vertices. 
Thus, there exists $z'\in N_{G[V_{k+2}]}(t)\setminus N_{G[V_{k+2}]}(s)$.

For each $sw\in E_{k+2}$, since $s\notin S_{sw}$ and $S_{sw}$ covers every edge of $G[V_{k+2}]-sw$, we have $N_{G[V_{k+2}]}(s)\setminus\{w\}\subseteq S_{sw}$. 
Since $tz'\in E_{k+2}$, we have $t\in S_{sw}$ or $z'\in S_{sw}$. 
Since $t,z'\notin N_{G[V_{k+2}]}(s)$, we have $|S_{sw}|\ge d_{G[V_{k+2}]}(s)$. 
By the same argument, $|S_e|\ge d_{G[V_{k+2}]}(t)$ for every edge $e\in E_{k+2}$ incident with $t$.

Since $st\notin E_{k+2}$, every edge of $G[V_{k+2}]$ is incident with exactly one of $s$ and $t$.
Summing over all edges of $E_{k+2}$, we have  $\sum_{e\in E_{k+2}}|S_e|\ge d_{G[V_{k+2}]}(s)^2+d_{G[V_{k+2}]}(t)^2\ge \frac{|E_{k+2}|^2}{2}\ge 3|E_{k+2}|-3$, 
where the last inequality follows from $|E_{k+2}|\ge5$. 
Thus, (ii) follows.   
\end{proof}

\begin{proof}[Proof of Theorem \ref{thm:main}]
Let $G:=G[A,B]$. 
Since the sharpness follows from Lemma \ref{lem:sharp}, it remains to prove that $|E(G)|\le(2k+1)(|A|-k)$. 
If $V_{k+2}=\emptyset$, then $|E(G)|=\frac{(k+1)|V(G)|}{2}=(k+1)|A|=(2k+1)|A|-k|A|\le
(2k+1)|A|-k(2k+1)=(2k+1)(|A|-k)$ as $|A|\ge3k$. 
Now we consider the case when $V_{k+2}\neq\emptyset$. 
Suppose to the contrary that $|E(G)|\ge(2k+1)(|A|-k)+1$. 
Then by Lemma \ref{lem:edge_count} (i) and $|V_{k+1}|=2|A|-|V_{k+2}|$, we have
\begin{equation}\label{eq:main_lowE1}
|E_{k+1}|\le |A|-(k+1)|V_{k+2}|+|E_{k+2}|+k(2k+1)-1.
\end{equation}
Combining Lemma \ref{lem:edge_count} (ii) and \eqref{eq:main_lowE1}, we have
\begin{equation}\label{eq:main_upE2}
|E_{k+2}|\ge k|V_{k+2}|-k(2k+1)+1+\min\{k,\nu(G[V_{k+2}])\}.
\end{equation} 
Since $V_{k+2}\neq\emptyset$, combining \eqref{eq:forest} and \eqref{eq:main_upE2}, we have $(k-1)|V_{k+2}|+\min\{k,\nu(G[V_{k+2}])\}\le (k-1)(2k+3)+1$ and so, $|V_{k+2}|\le 2k+3+\frac{1}{k-1}<2k+4$. 
Since $|V_{k+2}|$ is an integer,
\begin{equation}\label{eq:main_upV2}
|V_{k+2}|\le2k+3.
\end{equation}
Combining the lower bound on $2|E_{k+1}|$ in \eqref{eq:Se} and twice the upper bound on $|E_{k+1}|$ in
\eqref{eq:main_lowE1}, we have 
\begin{equation}\label{eq:main_E2Se}
\begin{aligned} 
(|A|-k-1-|V_{k+2}|)|E_{k+2}|+\sum_{e\in E_{k+2}}|S_e|\le &~(|V_{k+2}|-2k)|A|-2|A_{k+2}||B_{k+2}|\\
&~-(k+1)|V_{k+2}|+2k(2k+1)-2.
\end{aligned}
\end{equation}
If $|E_{k+2}|\ge2$ and $\nu(G[V_{k+2}])\ge2$, then $\sum_{e\in E_{k+2}}|S_e|\ge2|E_{k+2}|-2$ by \eqref{eq:sum_Se_matching}. 
If $|E_{k+2}|\ge2$ and $\nu(G[V_{k+2}])\le 1$, then $\nu(G[V_{k+2}])=1$ and $G[V_{k+2}]$ consists of a star with $|E_{k+2}|$ edges and possibly some isolated vertices. 
By \eqref{eq:sum_Se},
$\sum_{e\in E_{k+2}}|S_e|\ge |E_{k+2}|(|E_{k+2}|-1)\ge 2|E_{k+2}|-2$. 
If $|E_{k+2}|\le1$, then $|E_{k+2}|-1<0$ and $2|E_{k+2}|-2<0$.
If $|E_{k+2}|\ge2$, then $2|E_{k+2}|-2\ge|E_{k+2}|-1$.
Therefore, in all cases,
\begin{equation}\label{eq:main_sumSe}
\sum_{e\in E_{k+2}}|S_e|\ge\max\{|E_{k+2}|-1,\ 2|E_{k+2}|-2\}.
\end{equation}

Adjust notation if necessary, so that $|A_{k+2}|\le|B_{k+2}|$. 
\begin{claim}
    $k\le |A_{k+2}|\le |B_{k+2}|$.  
\end{claim}  
\begin{proof} 
Suppose to the contrary that $|A_{k+2}|\le k-1$.
Since every vertex of $B_{k+1}$ has degree $k+1$
and at most $|A_{k+2}|$ neighbors in $A_{k+2}$, we have
\begin{equation}\label{eq:main_lowE1-2}
|E_{k+1}|\ge(k+1-|A_{k+2}|)(|A|-|B_{k+2}|).
\end{equation}
Combining \eqref{eq:main_lowE1-2} and
\eqref{eq:main_lowE1}, we have
\begin{equation}\label{eq:main_lowE2}
\begin{aligned}
|E_{k+2}|\ge{}&
(k-|A_{k+2}|)|A|+|A_{k+2}||B_{k+2}|\\
&+(k+1)|A_{k+2}|-k(2k+1)+1.
\end{aligned}
\end{equation}
Since $G[V_{k+2}]$ is simple and bipartite,
we have $|E_{k+2}|\le|A_{k+2}||B_{k+2}|$.
By \eqref{eq:main_lowE2}, $|A|\ge3k$ and
$k-|A_{k+2}|>0$, we have
$(2k-1)|A_{k+2}|\ge k(k-1)+1>2k-1$,
where the last inequality follows from $k\ge3$.
Thus, $2\le|A_{k+2}|\le k-1$.

Combining \eqref{eq:main_lowE2} and \eqref{eq:forest},
and using $|V_{k+2}|=|A_{k+2}|+|B_{k+2}|$,
$|A|\ge3k$ and $k-|A_{k+2}|>0$, we have
\begin{equation}\label{eq:main_A2B2}
\begin{aligned}
(|A_{k+2}|-1)(2k-|B_{k+2}|)
&\ge(k-|A_{k+2}|)(|A|-3k)+(k-1)(k-2)\\
&\ge(k-1)(k-2).
\end{aligned}
\end{equation}
Since $1\le|A_{k+2}|-1\le k-2$, we have
$|B_{k+2}|\le k+1$ by \eqref{eq:main_A2B2}.
Hence, $|V_{k+2}|\le2k$.  
Since $|A|\ge 3k$ and $|V_{k+2}|\le 2k$, we have
$|A|-k-|V_{k+2}|\ge 0$. By \eqref{eq:main_sumSe},
$\sum_{e\in E_{k+2}}|S_e|\ge |E_{k+2}|-1$. Hence
\begin{equation}\label{eq:main_lowE2Se}
\begin{aligned}
&(|A|-k-1-|V_{k+2}|)|E_{k+2}|+\sum_{e\in E_{k+2}}|S_e| \\ 
\ge &~(|A|-k-1-|V_{k+2}|)|E_{k+2}|+|E_{k+2}|-1\\
=&~(|A|-k-|V_{k+2}|)|E_{k+2}|-1
\ge -1.
\end{aligned}
\end{equation}
Since $|V_{k+2}|-2k\le 0$ and $|A|\ge 3k$, we have
$(|V_{k+2}|-2k)|A|\le 3k(|V_{k+2}|-2k)$. 
Using \eqref{eq:main_E2Se} and 
$(|V_{k+2}|-2k)|A|\le 3k(|V_{k+2}|-2k)$, we have 
\begin{equation}\label{eq:main_upE2Se}
{\begin{aligned}
&(|A|-k-1-|V_{k+2}|)|E_{k+2}|+\sum_{e\in E_{k+2}}|S_e|\\
\le&~ (|V_{k+2}|-2k)|A|-2|A_{k+2}||B_{k+2}|-(k+1)|V_{k+2}|+2k(2k+1)-2\\
\le&~ 3k(|V_{k+2}|-2k)-2|A_{k+2}||B_{k+2}|-(k+1)|V_{k+2}|+2k(2k+1)-2\\
=&~(2k-1-2|A_{k+2}|)(|B_{k+2}|-k+1)+|A_{k+2}|-k-1,
\end{aligned}}
\end{equation} 
where the last equality follows from $|V_{k+2}|=|A_{k+2}|+|B_{k+2}|$. 
Combining \eqref{eq:main_lowE2Se} and
\eqref{eq:main_upE2Se}, since $|A_{k+2}|\le k-1$, we have 
\begin{equation*}
    (2k-1-2|A_{k+2}|)(|B_{k+2}|-k+1)\ge k-|A_{k+2}|\ge 1
\end{equation*} 
As $2k-1-2|A_{k+2}|\ge 1$, we have $|B_{k+2}|-k+1>0$, and hence $|B_{k+2}|\ge k$. 
Then $2k-|B_{k+2}|\le k$ and $|A_{k+2}|-1\ge 1$ and so, by \eqref{eq:main_A2B2}, we have $k(|A_{k+2}|-1)\ge (|A_{k+2}|-1)(2k-|B_{k+2}|)\ge (k-1)(k-2)=k(k-3)+2$.  
Thus, $|A_{k+2}|>k-2$. Since $|A_{k+2}|$ is an integer at most $k-1$, we have $|A_{k+2}|=k-1$.

Substituting $|A_{k+2}|=k-1$ into
\eqref{eq:main_lowE2}, and using $|A|\ge3k$
and $|B_{k+2}|\ge k$, we have
\begin{equation}\label{eq:main_lowE2-2}
|E_{k+2}|
\ge|A|+(k-1)|B_{k+2}|-k(k+1)
\ge k.
\end{equation}
Substituting $|A_{k+2}|=k-1$ into
\eqref{eq:main_upE2Se}, and using
\eqref{eq:main_sumSe},
\eqref{eq:main_lowE2-2},
$|A|\ge3k$, $|V_{k+2}|\le2k$ and $|B_{k+2}|\le k+1$,
we have
$1\le|E_{k+2}|-2
\le(|A|-k-1-|V_{k+2}|)|E_{k+2}|
+\sum_{e\in E_{k+2}}|S_e|
\le|B_{k+2}|-k-1\le0$,
a contradiction.
Therefore, the claim holds.
\end{proof}

If $|V_{k+2}|\ge2k+1$, then $|E_{k+2}|\ge1$
by \eqref{eq:main_upE2}.
Hence, $\nu(G[V_{k+2}])\ge1$, and
\eqref{eq:main_upE2} implies
\begin{equation}\label{eq:main_lowE2-3}
|E_{k+2}|\ge k(|V_{k+2}|-2k-1)+2.
\end{equation}

Suppose that $|V_{k+2}|=2k+3$. 
By \eqref{eq:main_lowE2-3} and \eqref{eq:forest}, we have $|E_{k+2}|=2k+2$.
Substituting $|V_{k+2}|=2k+3$ and $|E_{k+2}|=2k+2$
into \eqref{eq:main_upE2}, we have
$\min\{k,\nu(G[V_{k+2}])\}\le1$.
Since $k\ge3$ and $\nu(G[V_{k+2}])\ge1$, we have $\nu(G[V_{k+2}])=1$. 
Then $G[V_{k+2}]$ consists of a star with $2k+2$ edges and possibly some isolated vertices.
Since $|V_{k+2}|=2k+3$, we have $G[V_{k+2}]\cong K_{1,2k+2}$, contradicting the Claim. 
Thus, $|V_{k+2}|\neq2k+3$.
By \eqref{eq:main_upV2} and the Claim, we have $|V_{k+2}|\in\{2k,2k+1,2k+2\}$. 

If $|V_{k+2}|=2k$, then $|E_{k+2}|\ge 0=|V_{k+2}|-2k$. 
If $|V_{k+2}|\ge2k+1$, then $|E_{k+2}|\ge |V_{k+2}|-2k$ by \eqref{eq:main_lowE2-3}.
Therefore,
\begin{equation}\label{eq:main_lowE2-4}
|E_{k+2}|\ge|V_{k+2}|-2k.
\end{equation}

If $|V_{k+2}|=2k+2$, then $|E_{k+2}|\ge k+2\ge5$ by \eqref{eq:main_lowE2-3}.
Since $|V_{k+2}|\in\{2k,2k+1,2k+2\}$, using \eqref{eq:main_sumSe} for $|V_{k+2}|\in\{2k,2k+1\}$ and
Lemma \ref{lem:lowE1} (ii)  for $|V_{k+2}|=2k+2$, we have
\begin{equation}\label{eq:main_low_sumSe}
\sum_{e\in E_{k+2}}|S_e|\ge (|V_{k+2}|-2k+1)(|E_{k+2}|-1).
\end{equation}

By the Claim, we have $(|A_{k+2}|-k)(|B_{k+2}|-k)\ge0$.
Since $|A_{k+2}|+|B_{k+2}|=|V_{k+2}|$, we have 
\begin{equation}\label{eq:main_A2B2-2}
|A_{k+2}||B_{k+2}|\ge k(|V_{k+2}|-k).
\end{equation}
By \eqref{eq:main_low_sumSe}, \eqref{eq:main_E2Se} and \eqref{eq:main_A2B2-2}, we have
\begin{equation}\label{eq:main_lowE2-5}
\begin{aligned}
(|A|-3k)|E_{k+2}|-|V_{k+2}|+2k-1
&\le (|A|-k-1-|V_{k+2}|)|E_{k+2}|
+\sum_{e\in E_{k+2}}|S_e|\\
&\le (|V_{k+2}|-2k)(|A|-3k-1)-2.
\end{aligned}
\end{equation}
By \eqref{eq:main_lowE2-5}, we have $(|A|-3k)(|E_{k+2}|-|V_{k+2}|+2k)\le-1$.
However, $|A|\ge 3k$ and $|E_{k+2}|-|V_{k+2}|+2k\ge0$ by \eqref{eq:main_lowE2-4}, a contradiction.
Therefore, $|E(G)|\le(2k+1)(|A|-k)=(2k+1)(|V(G)|-2k)/2$.
\end{proof}

\section{Sharpness of $|V(G)|\ge 6k$}\label{sec:6k}
\begin{figure}[!h]
\centering
\includegraphics[height=3.5cm]{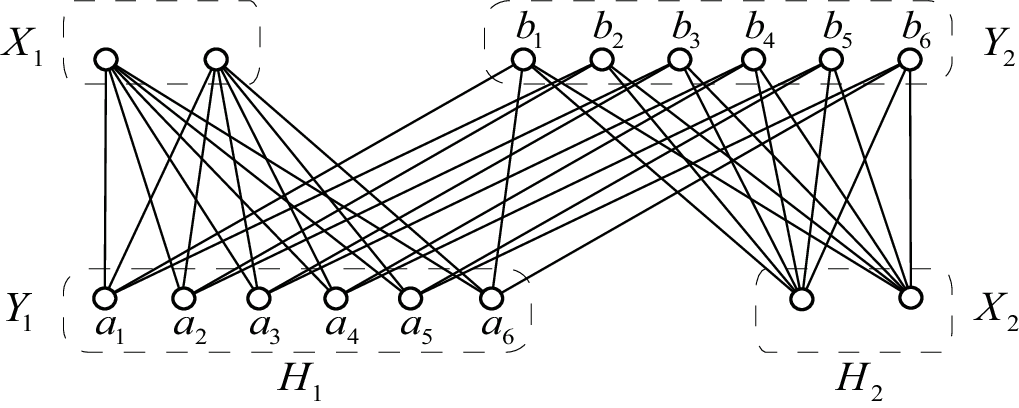}
\caption{The graph $H_{3,16}$.}
\label{fig:H_kn}
\end{figure}

For any integer $k\ge3$ and any even integer $n$ with $2k+2\le n\le6k-2$, let $t:=\lfloor(n-2k)/4\rfloor$. 
Let $H_1[X_1,Y_1]$ and $H_2[X_2,Y_2]$ be two disjoint copies of $K_{t,n/2-t}$ such that $|X_1|=|X_2|=t$, $Y_1:=\{a_1,\ldots,a_{n/2-t}\}$ and
$Y_2:=\{b_1,\ldots,b_{n/2-t}\}$.  
Let $H_{k,n}$ be the graph obtained from $H_1\cup H_2$ by joining each $a_i$ to the $k-t+1$ consecutive vertices $b_i,b_{i+1},\ldots,b_{i+k-t}$ of $Y_2$, for $1\le i\le n/2-t$, where the subscripts are taken modulo $n/2-t$ (see Figure \ref{fig:H_kn} for $H_{3,16}$). 
Then $H_{k,n}$ is bipartite with $A:=Y_1\cup X_2$ and
$B:=X_1\cup Y_2$ as its two partite sets.
In particular, if $t=0$, then $n=2k+2$ and
$H_{k,n}\cong K_{k+1,k+1}$. 

\begin{pro}\label{pro:6k}
For any integer $k\ge3$ and any even integer $n$ with $2k+2\le n\le6k-2$, $H_{k,n}$ is a minimal $k$-extendable bipartite graph of order $n$ with
$|E(H_{k,n})|>(2k+1)(n-2k)/2$.
\end{pro} 
\begin{proof}
Since $t=\lfloor(n-2k)/4\rfloor$ and $2k+2\le n\le6k-2$, we have 
\begin{equation}\label{eq:Hn}
0\le t\le k-1 ~\mbox{and}~ k+t\le n/2-t\le2k.
\end{equation}
Since $k-t\ge1$, we have $a_ib_i,a_ib_{i+1}\in E(H_{k,n})$ for $1\le i\le n/2-t$ and so, the cycle
$b_1a_1b_2a_2\ldots b_{n/2-t}a_{n/2-t}b_1$ is a spanning subgraph of $H_{k,n}[Y_1\cup Y_2]$. 
For $i=1,2$, every vertex of $X_i$ has a neighbor
in $Y_i$. Hence $H_{k,n}$ is connected. 
Moreover, $|A|=|B|=n/2$ and $|V(H_{k,n})|=n$. 
\begin{claim}
For every nonempty subset $Y'_1\subseteq Y_1$, $|N_{H_{k,n}}(Y'_1)\cap Y_2|\ge \min\{n/2-t,\,|Y'_1|+k-t\}$. 
\end{claim} 
\begin{proof}
If $Y_2\subseteq N_{H_{k,n}}(Y'_1)$, the claim follows.
Otherwise, by cyclically relabeling the vertices, we may assume that $b_{n/2-t}\notin N_{H_{k,n}}(Y'_1)$. 
Let $l$ be the largest index for which $a_l\in Y'_1$.
Since $b_{n/2-t}a_l\notin E(H_{k,n})$, we have $l+k-t<n/2-t$.
The set $N_{H_{k,n}}(Y'_1)\cap Y_2$ contains the $|Y'_1|$ vertices $b_j$ with $a_j\in Y'_1$, together with the $k-t$ further vertices $b_{l+1},\ldots,b_{l+k-t}$.
This proves the claim.
\end{proof}

To prove that $H_{k,n}$ is $k$-extendable, let $\emptyset\neq A'\subseteq A$.
If $A'\subseteq X_2$, then $N_{H_{k,n}}(A')=Y_2$, and so, $|Y_2|\ge k+t\ge k+|A'|$ by \eqref{eq:Hn}. 
If $A'\subseteq Y_1$, then $X_1\subseteq N_{H_{k,n}}(A')$, and $|N_{H_{k,n}}(A')|=t+|N_{H_{k,n}}(A')\cap Y_2|\ge \min\{|B|,\,|A'|+k\}$ by the Claim. 
If $A'\cap X_2\neq\emptyset$ and $A'\cap Y_1\neq\emptyset$, then $N_{H_{k,n}}(A')=X_1\cup Y_2=B$. 
Thus, $H_{k,n}$ is $k$-extendable by Theorem \ref{thm:char_k-ext-bi} (ii). 

For each $i\in\{1,2\}$, each vertex of $Y_i$ has $t$ neighbors in $X_i$ and $k-t+1$ neighbors in $Y_{3-i}$, and hence has degree $k+1$. 
For every $uv\in E(H_{k,n})$, at least one end vertex  of $uv$, say $u$, belongs to $Y_1\cup Y_2$. 
Hence $d_{H_{k,n}}(u)=k+1$ and $d_{H_{k,n}-uv}(u)=k$. 
So $H_{k,n}-uv$ is not $k$-extendable by Theorem \ref{thm:con_k-ext}.
Since $uv$ is arbitrary, $H_{k,n}$ is minimal $k$-extendable.

Since $|E(H_i)|=t(n/2-t)$ for each $i\in\{1,2\}$ and $|E_{H_{k,n}}[Y_1,Y_2]|=(k-t+1)(n/2-t)$, we have $|E(H_{k,n})|= |E(H_1)|+|E(H_2)|+|E_{H_{k,n}}[Y_1,Y_2]|=2t(n/2-t)+(k-t+1)(n/2-t)=(k+t+1)(n/2-t)$. 
Since $k-t\ge1$ and $2k+1-n/2+t\ge1$ by \eqref{eq:Hn}, we have $|E(H_{k,n})|-\frac{(2k+1)(n-2k)}2=(k-t)(2k+1-n/2+t)\ge1$. 
Thus, the result holds. 
\end{proof}

\end{document}